\documentclass[12pt]{article}

\usepackage{amsmath}

\usepackage{amsfonts}

\usepackage{amssymb}

\newtheorem{pro}{}

\newtheorem{lemma}[pro]{Lemma}

\makeatletter
\newfont{\footsc}{cmcsc10 at 8truept}
\newfont{\footbf}{cmbx10 at 8truept}
\newfont{\footrm}{cmr10 at 10truept}

\makeatother

\title{New upper bound for a class of vertex Folkman numbers}

\author{N.Kolev  \\%\thanks{Thanks to the editors of this wonderful journal!}\\
\small Department of Algebra\\[-0.8ex]
\small Faculty of Mathematics and Informatics\\[-0.8ex]
\small "St. Kl. Ohridski" University of Sofia\\[-0.8ex]
\small 5 J. Bourchier blvd, 1164 Sofia\\[-0.8ex]
BULGARIA \\[1ex]
N.Nenov  \\%\thanks{Thanks to the editors of this wonderful journal!}\\
\small Department of Algebra\\[-0.8ex]
\small Faculty of Mathematics and Informatics\\[-0.8ex]
\small "St. Kl. Ohridski" University of Sofia\\[-0.8ex]
\small 5 J. Bourchier blvd, 1164 Sofia\\[-0.8ex]
BULGARIA \\[-0.8ex]
\small \texttt{nenov@fmi.uni-sofia.bg}
}

\date{\small
Submitted: Jan 1, 2005;  Accepted: Jan 2, 2005; Published: Jan 3, 2005}

\begin{document}

\maketitle

\begin{abstract}
  Let $a_1$ ,\dots , $a_r$ be positive integers, $m=\sum_{i=1}^{r}
  (a_{i}-1)+1$ and $p= \max \{a_1, \ldots ,a_r\}$.  For a graph $G$
  the symbol $G\rightarrow \{a_1 , \ldots , a_r\}$ denotes that in
  every $r$-coloring of the vertices of $G$ there exists a
  monochromatic $a_i$-clique of color $i$ for some $i=1 , \ldots , r$.
  The vertex Folkman numbers $ F(a_1 , \dots , a_r ; m-1)=\min \{|
  V(G) | : G\rightarrow (a_1\dots a_r)$ and $K_{m-1} \not \subseteq G
  \}$ are considered. We prove that $F(a_1 , \ldots , a_r; m-1)\leq
  m+3p,$ $ p\geq 3 $. This inequality improves the bound for these
  numbers obtained in [6].
\end{abstract}

\section{Introduction}

We consider only finite, non-oriented graphs without loops and
multiple edges. We call a p-clique of the graph G a set of $p$
vertices, each two of which are adjacent.  The largest positive
integer p, such that the graph G contains a p-clique is denoted by
$cl(G)$ . In this paper we shall also use the following notations:

$V(G)-$ vertex set of the graph G;

$E(G)-$ edge set of the graph G;

$\bar{G}-$ the complement of G;

G[V], $V \subseteq V(G)-$ the subgraph of $G$ induced by $ V$;

$G-V - $ the subgraph induced by the set $ V(G)\backslash V;$

$N_G(v), v\in V(G)-$ the set of all vertices of G adjacent to v;

$K_n-$ the complete graph on $n$ vertices;

$C_n-$ simple cycle on n vertices;

$P_n- $ path on n vertices;

$\chi(G)-$ the chromatic number of $G;$

$\lceil x \rceil - $ the least positive integer greater or equal to
$x$.

Let $G_1$ and $G_2$ be two graphs without common vertices. We denote
by $G_1 + G_2$ the graph $G$ for which $V(G)=V(G_1) \cup V (G_2)$ and
$E(G)=E(G_1)\cup E(G_2) \cup E'$ , where $ E' = \{[x , y ] \mid x \in
V(G_1) , y \in V(G_2) \}$.

\vspace{5 mm}

\textbf{Definition} \textit{Let $ a_1 , \ldots , a_r $ be positive
  integers. We say that the $r$-coloring
   $$ V(G)= V_1 \cup \ldots \cup V_r \text{ , } V_i \cap V_j = \emptyset \text{ , } i \neq
   j$$ of the vertices of the graph $G$ is $(a_1 ,\ldots , a_r)$-free,
   if $V_i$ does not contain an $a_i$-clique for each $i \in \{ 1 ,
   \ldots , r \} $.}  \textit{ The symbol $G\rightarrow (a_1 , \ldots
   , a_r) $ means that there is not $(a_1 , \ldots , a_r)$ - free
   coloring of the vertices of $G.$} \vspace{5 mm}

 We consider for arbitrary natural numbers $a_1 , \ldots , a_r $ and
 $q$
$$ H(a_1,\ldots a_r ; q) = \{G :  G\rightarrow (a_1,\ldots
,a_r) \text { and } cl(G)<q \}$$

The vertex Folkman numbers are defined by the equalities
$$ F(a_1, \ldots , a_r ;q) = \min\{ | V(G)| : G\in H(a_1, \ldots ,
a_r; q)\}.$$ It is clear that $G\rightarrow (a_1,\ldots ,a_r)$ implies
$cl(G)\geq \max\{a_1,\ldots , a_r\}$. Folkman [3] proved that there
exists a graph $G$ such that $G \rightarrow (a_1, \ldots , a_r) $ and
$cl(G) = \max\{a_1,\ldots , a_r\}$. Therefore
\begin{equation}
  \label{1} F(a_1 ,\ldots , a_r ;q) \text{ exists
    if and only if }
  q> \mbox{max} \{ a_1,\ldots , a_r \}.
\end{equation}
These numbers are called vertex Folkman numbers. In [5] Luczak and
Urbansky defined for arbitrary positive integers $a_1 , \ldots , a_r$
the numbers
\begin{equation} \label{2} m= m(a_1 ,\ldots ,a_r) = \sum_{i=1}^{r}
  (a_i-1) +1 \text{ and } p=p(a_1,\ldots , a_r)=\max \{a_1 , \ldots ,
  a_r\}\text{.}
\end{equation}
Obviously $K_m\rightarrow(a_1,\ldots ,a_r)$ and $K_{m-1} \nrightarrow
(a_1,\ldots , a_r )$.  Therefore if $q\geq m+1$ then $F (a_1,\ldots
,a_r;q)= m $.

From (1) it follows that the number $F(a_1,\ldots ,a_r ;q)$ exists if
and only if $ m\geq p+1$.  Luczak and Urbansky [5] proved that
$F(a_1,\ldots ,a_r ;q)=m+p$. Later ,in [6], Luczak, Rucinsky and
Urbansky proved that $K_{m-p-1}+\bar{C}_{2p+1}$ is the only graph in
$H( a_1 ,\ldots , a_r;m)$ with $m+p$ vertices.

From (1) it follows that the number $F(a_1, \ldots ,a_r; m-1)$ exists
if and only if $m\geq p+2$. An overview of the results about the
numbers $F(a_1,\ldots ,a_r ;m-1) $ was given in [1]. Here we shall
note only the general bounds for the numbers $F(a_1,\ldots , a_r
;m-1)$. In [8] the following lower bound was proved $$F(a_1,\ldots
,a_r ;m-1) \geq m+p+2 \text{ , } p\geq2.$$ In the above equality an
equality occurs in the case when $\max \{a_1 ,\ldots ,a_r\}=2$ and
$m\geq 5$ (see [4,6,7]). For these reasons we shall further consider
only the numbers $F(a_1,\ldots ,a_r;m-1)$ when $\max\{a_1 , \ldots ,
a_r\}\geq3$.

In [6] Luczak, Rucinsky and Urbansky proved the following upper bound
for the numbers $F(a_1,\ldots ,a_r;m-1)$ :
$$F(a_1 , \ldots , a_r ;m-1)  \leq m+p^2 \text{ , for } m \geq 2p+2$$
In [6] they also announced without proof the following inequality :
$$F(a_1 , \ldots , a_r;m-1) \leq 3p^2+p-mp+2m-3  \text { , for } p+3 \leq m
\leq 2p+1.$$ In this paper we shall improve these bounds proving the
following \vspace{5 mm }

\textbf{ Main theorem}\textit{ Let $ a_1 ,\ldots , a_r $ be positive
  integers and $m$ and $p$ be defined by} (2)\textit{.  Let $m \geq p$
  and $p \geq 3.$ Then$$ F(a_1 , \ldots , a_r;m-1) \leq m+3p. $$}

     \textbf{Remark} \textit{This bound is exact for the numbers $F
       (2,2,3;4)$ and $F(3,3;4)$ because $$F(2,2,3;4)=14 \text{ (see
         [2]) }$$ and
    $$F(3,
    3;4)=14 \text{ (see [9]).}$$}

\section{Main construction}
We consider the cycle $C_{2p+1}$. We assume that $$ V( C_{2p+1})= \{
v_1 , \ldots , v_{2p+1}\}
$$ and $$ E( C _{2p+1}) = \{ [v_i , v_{i+1}], i=1 ,  \ldots , 2p\} \cup
\{ v_1 , v_{2p+1}\}.$$ Let $ \sigma$ denote the cyclic automorphism of
$C_{2p+1}$, i.e.  $\sigma(v_i)=v_{i+1}$ for $i=1 , \ldots , 2p,$
$\sigma(v_{2p+1})=v_1.$Using this automorphism and the set $M_1=
V(C_{2p+1}) \backslash \{v_1 , v_{2p-1} , v_{2p-2}\}$ we define $M_i =
\sigma^{i-1}(M_1) $ for $i=1 , \ldots , 2p+1$.  Let $ \Gamma_p$ denote
the extension of the graph $\bar{C}_{2p+1}$ obtained by adding the new
pairwise independent vertices $ u_1, \ldots , u_{2p+1} $ such that
\begin{equation}
  \label{3} N_{\Gamma_p} (u_i)=M_i \text{ for } i=1 , \ldots
  ,  2p+1  \text{ . }
\end{equation}
We easily see that $ cl( \bar{C}_{2p+1})=p$ .

Now we extend $\sigma$ to an automorphism of $\Gamma_p$ via the
equalities $\sigma(u_i)=u_{i+1}$ , for $i=1 , \ldots , 2p,$ and
$\sigma(u_{2p+1})=u_1$ Now it is clear that
\begin{equation}
  \label{4}     \sigma \text{ is an isomorphism of }
  \Gamma_p \text{.}
\end{equation}
The graph $\Gamma_p$ was defined for the first time in [ 8]. In [8]
it is also proved that $\Gamma_p \rightarrow (3,p)$ for $p \geq 3$.
For the proof of the main theorem we shall also use the following
generalisation of this fact.

\vspace{5 mm}

\textbf{Theorem 1.} \textit{ Let $p\geq 3$ be a positive integer}
\textit{and $m= p+2$. Then for arbitrary positive} \textit{ integers
  $a_1, \ldots , a_r$ ($r$ is not fixed)} \textit{such that
                    $$  m= 1+ \sum_{i}^{r} (a_i - 1)$$}

                  \textit{and $\max\{a_1 , \ldots , a_r\} \leq p $ we
                    have $$\Gamma _p \rightarrow (a_1 , \ldots a_r).$$
                  }

                  \section{Auxiliiary results}

                  The next proposition is well known and easy to
                  prove.

                  \vspace{5 mm }

                  \textbf{ Proposition 1} \textit{Let $a_1 ,\ldots ,
                    a_r $ be positive integers and $ n = a_1+ \ldots
                    +a_r$. Then} \textit{$$ \lceil \frac{a_1}{2}
                    \rceil + \ldots +\lceil \frac{a_r}{2} \rceil \geq
                    \lceil \frac{n}{2} \rceil .$$ If $n$ is even than
                    this inequality is strict unless all the numbers}
                  \textit{$ a_1 , \ldots , a_r $ are even. If $n$ is
                    odd then this} \textit{inequality is strict unless
                    exactly one of the numbers $a_1, \ldots , a_r$ is
                    odd.}

                  \vspace{5 mm} Let $ P_k $ be the simple path on $k$
                  vertices. Let us assume that
 $$V(P_k) = \{ v_1 , \ldots , v_k \}$$ and $$E(P_k)=\{ [ v_i , v_{i+1}], i=1 , \ldots , k-1\} .$$

 We shall need the following obvious facts for the complementary graph
 $ {\bar{P}_k}$ of the graph $P_k$:

 \begin{equation}
   \label{5}    cl(\bar{P}_k)= \lceil
   \frac{k}{2} \rceil
 \end{equation}
 \begin{equation}
   \label{6}     cl(\bar{P}_{2k}-v) = cl(\bar{P}_{2k}) \text{, for each } v\in V(
   \bar{P}_{2k})
 \end{equation}
 \begin{equation}
   \label{7}
   cl(\bar{P}_{2k}-\{ v_{2k-2} , v_{2k-1} \}) = cl(\bar{P}_{2k+1})
   \text{ for }
   k\geq 2
 \end{equation}
 \begin{equation}
   \label{8} cl(\bar{P}_{2k+1}-v_{2i}) = cl( \bar{P}_{2k+1})\text{ ,
   }i=1 ,
   \ldots , k \text{ , } k\geq 1.
 \end{equation}
 The proof of Theorem 1 is based upon three lemmas.

 \begin{lemma}

   Let $ V\subset V(C_{2p+1})$ and $ |V| = n< 2p+1.$ Let $ G=
   \bar{C}_{2p+1}[V] $ and let $ G_1 , \ldots , G_s$ be the connected
   components of the graph $\bar{G} = C_{2p+1}[V].$

   Then
   \begin{equation}
     \label{9}  cl(G)\geq  \lceil \frac{n}{2 } \rceil \text{.}
   \end{equation}
   If $n$ is even, then (9) is strict unless all $ |V(G_i)|$ for $i=1
   , \ldots , s$ are even.  If $n$ is odd, then (9) is strict unless
   exactly one of the numbers $|V(G_i)|$ is odd.
 \end{lemma}

 \textbf{ Proof} Let us observe that
 \begin{equation}
   \label{10}       G = \bar{G}_1 +\ldots +\bar{G}_s .
 \end{equation}

 Since $V \neq V( C_{2p+1})$ each of the graphs $G_i$ is a path.From
 (10) and (5) it follows that
     $$cl(G)= \sum_{i=1}^{s}\lceil \frac{n_i}{2} \rceil ,$$
     where $n_i= |V(G_i)|,$ $i=1 , \ldots , s .$ From this inequality
     and Proposition 1 we obtain the inequality (9) . From Proposition
     1 it also follows that if $n$ is even then there is equality in
     (9) if and only if the numbers $n_1 , \ldots , n_s$ are even ,
     and if $n$ is odd then we have equality in (9) if and only if
     exactly one of the numbers $n_1, \ldots , n_s$ is odd.

     \vspace{5 mm} \textbf{Corollary 1} \textit{It is true that
       $cl(\Gamma_p) =p.$ } \vspace{ 5 mm}

     \textbf{Proof} It is obvious that $cl(\bar{C}_{2p+1})=p$ and
     hence $cl(\Gamma_p)\geq p.$ Let us denote an arbitrary maximal
     clique of $\Gamma_p$ by $Q.$ Let us assume that $|Q|> p.$ Then
     $Q$ must contain a vertex $u_i$ for some $i=1 , \ldots , 2p+1.$
     As the vertices $u_i$ are pairwise independent $Q$ must contain
     at most one of them. Since $\sigma$ is an automorphism of
     $\Gamma_p$ (see (4)) and $u_i = \sigma^{i-1}(u_1),$ we may assume
     that $Q$ contains $u_1.$ Let us assign the subgraph of $\Gamma_p$
     induced by $ N_{\Gamma_p(u_1)}= M_1 $ by $H.$ The connected
     components of $H$ are $ \{ v_2 , v_3 , \ldots , v_{2p-3}\} $ and
     $ \{v_{2p} , v_{2p+1} \} $ and the both of them contain an even
     number of vertices. Using Lemma 1 we have $cl(H)=p-1.$ Hence
     $|Q|=p$ and this contradicts the assumption .

     \vspace { 5 mm }

     The next two lemmas follow directly from (10) , (6) , (7) , and
     (8) and need no proof.

     \begin{lemma} Let $ V\subsetneq V(C_{2p+1})$ and $G=
       \bar{C}_{2p+1}[V] .$ Let $P_k= \{v_1 , v_2 , \ldots , v_k\}$ be
       a connected component of the graph $ \bar{G} = C_{2p+1}[V].$
       Then

       (a) if $k=2s$ then
        $$cl(G-v_i)=cl(G),  i=1 , \ldots , 2s$$
        and
        $$cl(G- \{v_{2s-2} , v_{2s-1} \}) = cl(G).$$
        (b) if $k=2s+1$ then $$cl(G-v_{2i}) = cl(G) , i=1 , \ldots ,
        s$$
      \end{lemma}
      \begin{lemma} Let $V\subseteq V(C_{2p+1}) $ and $ \bar{C}_{2p+1}
        =G.$ Let $ P_{2k}= \{v_1 , \ldots , v_{2k} \} $ and $ P_s = \{
        w_1 , \ldots , w_s \}$ be two connected components of the
        graph $ \bar{G} = C_{2p+1} [V] .$ Then

        (a) if $s=2t$ then $$cl(G-\{v_i , w_j\}) = cl(G), $$ for $i=1
        , \ldots , 2k, $ $j= 1 , \ldots , s $ and
    $$cl(G-\{ v_{2k-2}, v_{2k-1}, w_j\})= cl(G),$$
    for $j=1 , \ldots , s.$

    (b) If $s=2t+1$ then $$cl(G- \{ v_{2k-2}, v_{2k-1} , w_{2i} \}) =
    cl(G) \text { , for } i=1 , \ldots , t.$$
  \end{lemma}

  \section {Proof of theorem 1}

  We shall prove Theorem 1 by induction on $r.$ As $m=\sum_{i=1}^{r}
  (a_{i}-1) +1 = p+2$ and $\max\{ a_1 , \ldots , a_r\} \leq p$ we have
  $r\geq2.$ Therefore the base of the induction is $r=2.$ We warn the
  reader that the proof of the inductive base is much more involved
  then the proof of the inductive step.  Let $r=2$ and $(a_1-1)
  +(a_2-1) +1=p+2$ and $\max\{a_1,a_2\}\leq p.$ Then we have
  \begin{equation}
    \label{11}  a_1+a_2=p+3.
  \end{equation}

  Since $p\geq3$ and $\max \{ a_1 , a_2 \} \leq p$ we have that
  \begin{equation}
    \label{12} a_i \geq 3 \text{ , } i=1,2 \text{.}
  \end{equation}
  We must prove that $\Gamma_p \rightarrow (a_1,a_2).$ Assume the
  opposite and let $V(\Gamma_p)= V_1\cup V_2$ be a $(a_1, a_2 )$-free
  coloring of $V(\Gamma_p).$ Define the sets $$V'_i=V_i \cap
  V(\bar{C}_{2p+1})\text{ , } i=1,2$$ and the graphs
               $$G_i= \bar{C}_{2p+1}[V'_i] \text{ , }
               i=1,2.$$ By assumption $V_i $ does not contain an
               $a_i$-clique and hence $V'_i$ does not contain an
               $a_i$-clique, too. Therefore from Lemma 1 we have
               $V'_i\leq 2a_i -2,$ i=1,2. From these inequalities and
               the equality
$$| V'_1 | + | V'_2 | = 2p+1 =2a_1+2a_2-5$$ (as
$p=a_1+a_2-3-$ see (11) ) we have two possibilities:
$$| V'_1 |=2a_1-2 \text{ , } |V'_2|=2a_2-3$$
or
$$| V'_1|=2a_1-3 \text{ ,} |V'_2|=2a_2-2.$$
Without loss of generality we assume that
\begin{equation}
  \label{13} | V'_1 |=2a_1-2 \text{ , } |V'_2|=2a_2-3 \text{.}
\end{equation}
From (13) and Lemma 1 we obtain $cl(G_i)\geq a_i-1$ and by the
assumption that the coloring $V_1\cup V_2$ is $(a_1 , a_2)$-free we
have
\begin{equation}
  \label{14} cl(G_i) =a_i - 1 \text{ for } i=1,2\text{.}
\end{equation}
From (13) , (14) and Lemma 1 we conclude that

\begin{equation} \label{15} \vcenter{\vbox{\hsize=10cm
      \begin{center}\it
        The number of the vertices of each connected component of\,
        $\bar{G}_1$ is an even number;
      \end{center}}}
\end{equation}
and
\begin{equation}
  \label{16}
  \vcenter{\vbox{\hsize=10cm
      \begin{center}\it
        the number of the vertices of exactly one of the connected
        components of\, $\bar{G}_2$ is an odd number.
      \end{center}}}
\end{equation}
According to (15) there are two possible cases.

\textbf{Case 1}. Some connected component of $\bar{G}_1$ has more then
two vertices. Now from (15) it follows that this component has at
least four vertices. Taking into consideration (15) and (4) we may
assume that $ \{ v_1 , \ldots , v_{2s}\} ,$ $s\geq 2$ is a connected
component of $\bar{G}_1.$ Since $V'_1$ does not contain an
$a_1$-clique we have by Lemma 1 that $s< a_1.$ Therefore $2s+2\leq 2p$
and we can consider the vertex $u_{2s+2}.$

\textbf{Subcase 1.a}. Assume that $u_{2s+2}\in V_1.$ Let $v_{2s+2}\in
V'_2.$ We have from (3) that
\begin{equation}
  \label{17}  N_{\Gamma_p}(u_{2s+2}) \supseteq V'_1-
  \{v_{2s-2},v_{2s-1} \}.
\end{equation}
From (14) and Lemma 2 (a) we have that $ V'_1-\{v_{2s-2} , v_{2s-1}\}$
contains an $(a_1-1)$-clique $Q$ . From (17) it follows that $Q\cup
\{u_{2s+2} \}$ is an $a_1$-clique in $V_1$ which is a contradiction.

Now let $v_{2s+2}\in V'_1.$ From (3) we have
\begin{equation}
  \label{18}  N_{\Gamma_P}(u_{2s+2})\supseteq V'_1 -
  \{ v_{2s-2} , v_{2s-1} , v_{2s+2} \}.
\end{equation}
According to (15) we can apply Lemma 3 (a) for the connected component
$\{v_1 , \ldots , v_{2s} \} $ of $\bar{G}_1$ and the connected
component of $\bar{G}_1$ that contains $v _{2s+2}.$ We see from (14)
and Lemma 3(a) that $V_1'-\{v_{2s-2} , v_{2s-1} , v_{2s+2} \}$
contains an $(a_1-1)$-clique $Q$ of the graph $G_1.$ Now from (18) it
follows that $Q\cup \{u_{2s+2} \}$ is an $a_1$-clique in $V_1,$ which
is a contradiction.

\textbf{Subcase 1.b}. Assume that $u_{2s+2}\in V_2.$ If
$v_{2s+2}\notin V'_2$ then from (3) it follows
\begin{equation}
  \label{19}  N_{\Gamma_p}(u_{2s+2})\supseteq V'_2.
\end{equation}
As $V'_2$ contains an $(a_2-1)$-clique Q (see (14) ). From (19) it
follows that $Q\cup \{u_{2s+2} \}$ is an $a_2$-clique in $V_2, $ which
is a contradiction.

Let now $v_{2s+2}\in V'_2.$ In this situation we have from (3)
\begin{equation}
  \label{20}  N_{\Gamma_p}(u_{2s+2})\supseteq
  V'_2-\{ v_{2s+2} \}.
\end{equation}
We shall prove that
\begin{equation}
  \label{21}  V_2 - \{ v_{2s+2} \} \text{ contains an }
  (a_2-1)\text{-clique of } \Gamma_p.
\end{equation}
As $v_{2s}$ is the last vertex in the connected component of $G_1,$ we
have $v_{2s+1} \in V'_2.$ Let $L$ be the connected component of
$\bar{G}_2$ containing $v_{2s+2}.$ Now we have $L= \{ v_{2s+1},
v_{2s+2} , \ldots \} .$ Now (21) follows from Lemma 2(b) applied to
the component $L.$ From (20) and(21) it follows that $V_2$ contains aa
$a_2$-clique, which is a contradiction.

\textbf{ Case 2}. Let some connected component of $\bar{G}_1$ have
exactly two vertices.

From (12) and (13) it follows that $\bar{G}_1$ has at least two
connected components. It is clear that $\bar{G}_2$ also has at least
two components. From (16) we have that the number of the vertices of
at least one of the components of $G_2$ is even. From these
considerations and (4) it follows that it is enough to consider the
situation when $ \{ v_1 , v_2 \}$ is a connected component of
$\bar{G}_1$ and $ \{v_3, \ldots , v_{2s} \} $ is a component of
$\bar{G}_2,$ and $ \{v_{2s+1} , v_{2s+2} \} $ is a component of
$\bar{G_1}.$ We shall consider two subcases.

\textbf{Subcase 2.a.} If $u_{2s+2} \in V_1.$

Let $s=2.$ We apply Lemma 3(a) to the components $\{ v_1 , v_2 \}$ and
$ \{ v_5, v_6 \}.$ From (14) we conclude that

\begin{equation}
  \label{22}  V'_1 - \{ v_2, v_6 \}  \text{
    contains an } (a_1 - 1)\text{-clique.}
\end{equation}

From (3) we have

\begin{equation}
  \label{23} N_{\Gamma_p}(u_6) \supseteq
  V'_1 -\{v_2,v_6\}.
\end{equation}
\vspace{3 mm}

Now (22) and (23) give that $V_1$ contains an $a_1$-clique.

Let $s\geq3.$ From (3) we have
\begin{equation}
  \label{24}  N_{\Gamma_p}(u_{2s+2}) \supseteq V'_1- \{ v_{2s+2} \}.
\end{equation}
According to Lemma 2(a) $V'_1- \{v_{2s+2} \}$ contains an
$(a_1-1)$-clique. Now using (24) we have that this $(a_1-1)$-clique
together with the vertex $u_{2s+2} $ gives an $a_1$-clique in $V_1.$
Subcase 2.a. is proved.

\textbf{Subcase 2.b}. Let $u_{2s+2}\in V_2.$

Let $s=2.$ From (3) we have $ N_{\Gamma_p}(u_6) \supseteq V'_2-\{ v_3
\}.$ According to Lemma 2(a) and (14) $V'_2- \{v_3\}$ contains an
$(a_2-1)$-clique. This clique together with $u_{2s+2}\in V_2$ gives an
$a_2$-clique in $V_2,$ which is a contradiction.

Let $s\geq3.$ Here from (3) we have $ N_{\Gamma_p}(u_{2s+2})\supseteq
V'_2- \{v_{2s-2} , v_{2s-1} \}.$ According to Lemma 2 (a) and (14) we
have that $V'_2 - \{v_{2s-2} , v_{2s-1} \}$ contains an
$(a_2-1)$-clique. This clique together with $u_{2s+2} \in V_2$ gives
an $a_2$-clique in $V_2,$ which is a contradiction.  This completes
the proof of case 2 and of the inductive base $r=2.$

Now we more easily handle the case $r\geq3.$ It is clear that
                $$G\rightarrow (a_1, \ldots , a_r) \Leftrightarrow
                G\rightarrow (a_{\varphi (1) } , \ldots , a_{ \varphi
                  (r)})$$ for any permutation $\varphi \in S_r.$ That
                is why we may assume that
                \begin{equation}
                  \label{25}  a_1 \leq \ldots \leq a_r \leq p.
                \end{equation}

                We shall prove that $a_1+a_2-1 \leq p.$ If $a_2\leq2$
                this is trivial : $a_1+a_2 -1 \leq 3 \leq p.$ Let
                $a_2\geq3.$ From (25) we have $a_i \geq3,$ $i=2,
                \ldots ,r.$ From these inequalities and the statement
                of the theorem
$$\sum_{i=1}^{r} (a_i-1)  +1=p+2$$
we have
 $$p+2 \geq 1+(a_2-1)+(a_1-1) + 2(r-2).$$
 From this inequality and $ r\geq 3$ it follows that $a_1 +a_2 -1 \leq
 p.$ Thus we can now use the inductive assumption and obtain
 \begin{equation}
   \label{26} \Gamma_p \rightarrow (a_1+a_2-1 , a_3 , \ldots , a_r).
 \end{equation}
 Consider an arbitrary $r$-coloring $V_1 \cup \ldots \cup V_r$ of
 $V(\Gamma_p).$ Let us assume that $V_i$ does not contain an
 $a_i$-clique for each $i=3 , \ldots , r.$ Then from (26) we have $V_1
 \cup V_2$ contains $(a_1 + a_2 -1)$-clique. Now from the pigeonhole
 principle it follows that either $V_1$ contains an $a_1$-clique or
 $V_2$ contains an $a_2$-clique. This completes the proof of Theorem
 1.

 \section{Proof of the main theorem}

 Let $m$ and $p$ be positive integers $p\geq3$ and $m\geq p+2.$ We
 shall first prove that for arbitrary positive integers $ a_1 , \ldots
 , a_r $ such that
$$ m= 1+ \sum_{i=1}^{r} (a_i - 1) $$

and $\max\{a_1 , \ldots , a_r\}\leq p$ we have
\begin{equation}
  K_{m-p-2} + \Gamma_p \rightarrow (a_1 , \ldots , a_r ).
\end{equation}

We shall prove (27) by induction on $t=m-p-2.$ As $m\geq p+2$ the base
is $t=0$ and it follows from Theorem 1. Assume now $t\geq 1.$ Then
obviously

    $$K_{m-p-2} + \Gamma_p = K_1 +( K_{m-p-3}+ \Gamma_p).$$

    Let $V(K_1)=\{ w \}.$ Consider an arbitrary $r$-coloring $V_1 \cup
    \ldots \cup V_r$ of $V( K_{m-p-2} + \Gamma_p).$ Let $w\in V_i$ and
    $V_j$ , $j\neq i$ does not contain a $a_j$-clique.

    In order to prove (27) we need to prove that $V_i$ contains an
    $a_i$-clique. If $a_i =1$ this is clear as $w\in V_i.$ Let $ a_i
    \geq 2.$ According to the inductive hypothesis we have
    \begin{equation}
      \label{28} K_{m-p-3}+ \Gamma_p \rightarrow (a_1 , \ldots
      ,a_{i-1},
      a_i -1 , a_{i+1} , \ldots , a_r).
    \end{equation}
    We consider the coloring
            $$V_1 \cup \ldots \cup V_{i-1} \cup \{V_i-w \}
            \cup \ldots \cup V_r$$ of $V(K_{m-p-3} + \Gamma_p).$ As
            $V_j,$ $j\neq i$ do not contain $a_j$-cliques, from (28)
            we have that $V_i-\{w\}$ contains an $(a_i-1)$-clique.
            This $(a_i-1)$-clique together with $w$ form an
            $a_i$-clique in $V_i.$ Thus (27) is proved.

            From Corollary 1 obviously follows that $cl(K_{m-p-2} +
            \Gamma_p) =m-2.$ From this and (27) we have $K_{m-p-2}+
            \Gamma_p \in H(a_1, \ldots , a_r;m-1).$ The number of the
            vertices of the graph $K_{m-p-2} + \Gamma_p$ is $ m+3p$
            therefore $F(a_1 , \ldots , a_r ; m-1) \leq m+3p.$

            The main theorem is proved.

\end{document}